\documentclass[12pt]{article}
\usepackage{amsmath,amsthm,amssymb,amscd}
\usepackage{latexsym}
\newtheorem{thm}{Theorem}[section]

 \newtheorem{prop}[thm]{Proposition}
 \theoremstyle{definition}
 
 \theoremstyle{remark}

 \numberwithin{equation}{section}

\title
{Three-orbifolds with positive scalar curvature}

\author{ Hong Huang}
\date{}
\begin{document}
\maketitle
\begin{abstract}
 We prove the following result: Let $(\mathcal{O},g_0)$ be a complete,
connected 3-orbifold with uniformly positive scalar curvature, with
bounded geometry,  and containing no bad 2-suborbifolds. Then there is a
finite collection $\mathcal{F}$ of spherical 3-orbifolds, such that
$\mathcal{O}$ is diffeomorphic to a (possibly infinite) orbifold
connected sum of copies of members in $\mathcal{F}$. This extends
work of Perelman and Bessi$\grave{e}$res-Besson-Maillot. The proof
uses Ricci flow with surgery on complete 3-orbifolds, and are along
the lines of the author's previous work on 4-orbifolds with positive isotropic curvature.


\end{abstract}
\maketitle


\section {Introduction}

In Perelman [P1, P2, P3] Perelman used Hamilton's Ricci flow to
attack Poincar$\acute{e}$ and geometrization conjectures; cf. also
[BBB$^+$], [CaZ], [KL1], [MT]  and [Z]. In particular he completed
the previous work of Gromov-Lawson (see for example [GL]) and Schoen-Yau [SY] on the
classification of compact, oriented 3-manifolds with positive scalar
curvature. In [BBM], Bessi$\grave{e}$res, Besson, and Maillot
extended Perelman's work to get a classification of open, oriented
3-manifolds  with uniformly positive scalar curvature and with
bounded geometry. In this short note we will go further to extend
the above work to the orbifold case. More precisely we will prove
the following

 \begin{thm} \label{thm 1.1}
  Let $(\mathcal{O},g_0)$ be a complete,
connected 3-orbifold with uniformly positive scalar curvature, with
bounded geometry,  and containing no bad 2-suborbifolds. Then there is a
finite collection $\mathcal{F}$ of spherical 3-orbifolds, such that
$\mathcal{O}$ is diffeomorphic to a (possibly infinite) orbifold
connected sum of copies of members in $\mathcal{F}$.
\end{thm}

 \noindent Note that the result in the compact, oriented case is implicit in [KL2]. 
 For an introduction to orbifold see [T]. Recall that an orbifold is good if it has a covering which is a manifold, otherwise it is bad.
 A Riemannian orbifold $(\mathcal{O}, g)$ is said to have
uniformly positive scalar curvature  if there exists a constant $c>0$ such that  the scalar curvature
 satisfies
$R\geq c$.
Also recall that a complete Riemannian orbifold $\mathcal{O}$ is
said to have bounded geometry if the sectional curvature is bounded
(in both sides) and the volume of any unit ball in  $\mathcal{O}$ is
uniformly bounded away from zero.  The notion of  a (possibly infinite) orbifold connected sum will be
given later in this section; cf. [Hu2].

By extending the construction in [GL] to the orbifold case and
utilizing a result in Dinkelbach-Leeb (see Proposition 2.2 in [DL]),
one can show that the converse of Theorem 1.1 is also true:
 Any 3-orbifold as in the conclusion of the theorem admits a complete metric
  with uniformly positive scalar curvature and with bounded
geometry.

 Now we explain the notion of (possibly infinite) orbifold connected sum (see [Hu2]),
 which extends several related  notions  in the literature (see for example [CTZ]  and [KL2]). Let $\mathcal{O}_i$ ($i=1,2$) be two $n$-orbifolds, and let
$D_i \subset \mathcal{O}_i$ be two suborbifolds-with boundary, both
diffeomorphic to some quotient orbifold $D^n//\Gamma$, where $D^n$
is the closed unit $n$-ball, and $\Gamma$ is a finite subgroup of
$O(n)$. Choose a diffeomorphism $f:
\partial D_1 \rightarrow \partial D_2$, and use it to glue together
$\mathcal{O}_1\setminus int (D_1 )$ and  $\mathcal{O}_2\setminus int
(D_2 )$. The result is called the orbifold connected sum of
$\mathcal{O}_1$ and $\mathcal{O}_2$ by gluing map $f$, and is
denoted by $\mathcal{O}_1 \sharp_f  \mathcal{O}_2 $. If $ D_i$
($i=1,2$) are disjoint suborbifolds-with boundary (both
diffeomorphic to some quotient orbifold $D^n//\Gamma$) in the same
$n$-orbifold $\mathcal{O}$, the result of similar process as above
is called the orbifold connected sum of $\mathcal{O}$  to itself,
and is denoted by $\mathcal{O}\sharp_f$.

  Given a collection $\mathcal{F}$ of $n$-orbifolds, we say an $n$-orbifold $\mathcal{O}$ is a (possibly infinite)
 orbifold connected sum of  members of $\mathcal{F}$  if there exist a  graph $G$ (in which we allow an edge to connect some vertex
to itself), a map $v\mapsto F_v$ which associates to each vertex of
$G$ a copy of some orbifold in $\mathcal{F}$, and a map $e \mapsto
f_e$ which associates to each edge of $G$ a self-diffeomorphism of
some $(n-1)$-dimensional spherical orbifold, such that if we do an
orbifold connected sum of $F_v$'s along each edge $e$ using the
gluing map $f_e$, we obtain an $n$-orbifold diffeomorphic to
$\mathcal{O}$.

 The proof of the theorem is along the lines of [Hu1, Hu2], which in turn are based on [P2] and
  which incorporate ideas from [BBB$^+$], [BBM], [CZ], [CTZ] and [H97].
 Meanwhile we will establish the short time existence of Ricci flow on complete orbifolds with bounded curvature, which extends Shi [S].
 This result was taken for granted of in the literature without a proof  (see for example [Hu2]).
   We will construct Ricci flow with surgery with initial data a  complete 3-orbifold with bounded geometry  and with no bad 2-suborbifolds,
 which extends [P2], [BBM] and [KL2].  In contrast to the oriented case considered in [P2], [BBM] and [KL2], new orbifold singularities may
 be introduced in the process of surgery in the nonoriented case.  Following [CTZ] and [Hu2], we deal with this situation using
 Hamilton's canonical parametrization for a neck region ([H97]), which also applies to the 3-dimensional case.

In Section 2 we  establish the short time existence of Ricci flow on complete orbifolds with bounded curvature, and describe the
canonical neighborhood structure of
 ancient $\kappa$-solutions on 3-orbifolds. In Section 3
we construct Ricci flow with surgery with initial data a  complete
3-orbifold with bounded geometry  and containing no bad 2-suborbifolds,
then Theorem 1.1 follows immediately.

\section{ Ancient $\kappa$-solutions on 3-orbifolds}

First we establish short time existence for Ricci flow on complete
orbifolds (not necessarily of 3-dimensional) with bounded curvature.
Let $(\mathcal{O},g_0)$ be a complete $n$-orbifold
 with $|Rm|\leq K$  for a constant $K$.
Consider the Ricci flow ([H82])

\begin{displaymath}\frac{\partial g}{\partial t}=-2Ric, \ \  g|_{t=0}=g_0. \ \ \ \
 \end{displaymath}
We will extend Shi's result in [S] to get the following

\begin{thm} \label{thm 2.1}  Let $(\mathcal{O},g_0)$ be a complete $n$-orbifold
 with $|Rm|\leq K$.  Then the Ricci flow with initial data $(\mathcal{O},g_0)$
 has a short time solution.
\end{thm}

{\bf Proof}  The proof of Shi [S] in the manifold case can be
adapted to the orbifold case with some minor modifications. We only
indicate the necessary changes. First one can write the complete
$n$-orbifold $\mathcal{O}$ as a union of a sequence of compact
$n$-suborbifolds with boundary $(n-1)$-suborbifolds. (This is
possible, as is seen by smoothing the distance function from a given
smooth point via heat equation, and using Sard theorem and preimage
theorem (see [BB]).) Also note that Stokes theorem holds in a
bounded domain (with boundary a suborbifold) in a Riemannian
orbifold. (In fact it holds for a slightly more general domain, see
[C]. ) In Lemma 3.1 (on p. 244) of [S], the dependence on the
injectivity radius of $\overline{D}$ can be replaced by that on the Sobolev
constant in some Sobolev inequality which holds in a bounded domain
(with boundary a suborbifold) in a Riemannian orbifold. (For Sobolev
inequalities on compact manifolds with boundary, one can see [A] and
[He]; for Sobolev inequalities  on closed orbifold, see [N] and [F].
The extension to the case of compact orbifolds with boundary is
routine.) In pages 260 and 286  of [S], one can pull back the solution to
(a suitable ball in) $T_{\tilde{x}_0}\widetilde{U}$  via
$exp_{x_0}\circ \pi_*$, where $(\widetilde{U}, G, \pi)$ is a
orbifold chart for some $U$ around $x_0$.
 \hfill{$\Box$

 \vspace *{0.4cm}

  By extending the proof of [Ko] to the orbifold case,
one also has uniqueness for the Ricci flow  (with a given initial
data) in the category of bounded curvature solutions  (even in a
slightly larger category).

From now on we will restrict to 3-dimensional orbifolds.  The notion of
ancient $\kappa$-solutions can be extended obviously from manifolds
to orbifolds. (See [KL2].) We will investigate  their structure.
First we define necks and caps. (See also [KL2].)  Let $\Gamma$ be a
finite subgroup of $O(3)$. A (topological) neck is an orbifold which
is diffeomorphic to $\mathbb{S}^2//\Gamma \times \mathbb{R}$, where
$\mathbb{S}^2//\Gamma $ denotes the quotient orbifold. To define
caps, let $\Gamma'$ be a finite subgroup of $O(3)$ such that
$\mathbb{S}^2//\Gamma'$ admits an isometric involution $\sigma$
(i.e. a nontrivial isometry $\sigma$ with $\sigma^2=1$), and
consider the quotient $(\mathbb{S}^2//\Gamma' \times
\mathbb{R})//\{1,\hat{\sigma}\}$, where $\hat{\sigma}$ is the
reflection on the orbifold $\mathbb{S}^2//\Gamma' \times \mathbb{R}$
defined by $\hat{\sigma}(x,s)=(\sigma(x),-s)$ for $x\in
\mathbb{S}^2//\Gamma' $ and $s\in \mathbb{R}$. We denote this
orbifold by $\mathbb{S}^2//\Gamma' \times_{Z_2}  \mathbb{R}$. (By
the way, note that $\Gamma'$ and $\hat{\sigma}$ may act
isometrically on $\mathbb{S}^3$ in a natural way.) We define a
(topological) cap to be an orbifold diffeomorphic to either
$\mathbb{S}^2//\Gamma' \times_{Z_2}  \mathbb{R}$ as above, or
$\mathbb{R}^3//\Gamma$, where $\Gamma$ is a finite subgroup of
$O(3)$.

It turns  out that any neck or cap can be written as  an
infinite orbifold connected sum of spherical 3-orbifolds.
Let $\Gamma$ be a finite subgroup of $O(3)$, then $\Gamma$ acts on
$\mathbb{S}^3$ by suspension.
 Using a somewhat ambiguous notation,  $\mathbb{S}^2//\Gamma \times \mathbb{R}^1 \approx
\cdot\cdot\cdot \sharp \mathbb{S}^3//\Gamma \sharp
\mathbb{S}^3//\Gamma \sharp \cdot\cdot\cdot$, $\mathbb{R}^3//\Gamma \approx
\mathbb{S}^3
//\Gamma \sharp \mathbb{S}^3 //\Gamma \sharp \cdot\cdot\cdot$, and
$\mathbb{S}^2//\Gamma' \times_{Z_2} \mathbb{R} \approx
\mathbb{S}^3//\{\Gamma',\hat{\sigma}\}\sharp \mathbb{S}^3//\Gamma'
\sharp \mathbb{S}^3//\Gamma' \sharp \cdot\cdot\cdot$. Note that the
orbifold connected sums appeared in these  three examples are
actually the operation of performing 0-surgery as defined in [KL2]
(which can be extended to the non-oriented case), and we have
omitted the $f$'s in the notation.  Also note that for a
diffeomorphism $f: \mathbb{S}^2//\Gamma \rightarrow
\mathbb{S}^2//\Gamma$  the mapping torus $\mathbb{S}^2//\Gamma
\times_f \mathbb{S}^1 \approx \mathbb{S}^3//\Gamma \sharp_f$.  By a
result in [DL] (see Proposition 2.2 there), $f$ is isotopic to an
isometry, so $\mathbb{S}^2//\Gamma \times_f \mathbb{S}^1 $ admits a
$\mathbb{S}^2 \times \mathbb{R}$-geometry.

\vspace*{0.4cm}

The following proposition is an analog of  Proposition 2.4 in [Hu2].

\begin{prop} \label{prop 2.2} \ \ There exist a universal positive constant $\eta$ and  a universal positive function $\omega: [0,\infty)\rightarrow (0,\infty)$
such that for any 3-dimensional ancient $\kappa$-solution
$(\mathcal{O},g(t))$, we have

\begin{displaymath}
(i) \hspace{4mm}  R(x,t)\leq R(y,t)\omega(R(y,t)d_t(x,y)^2)
\end{displaymath}
\noindent for any $x,y\in \mathcal{O}$, $t\in (-\infty,0]$, and

\begin{displaymath}(ii) \hspace{4mm}  |\nabla R|(x,t)\leq \eta R^{\frac{3}{2}}(x,t), \hspace{4mm}  |\frac{\partial R}{\partial
t}|(x,t)< \eta R^2(x,t)
\end{displaymath}
\noindent for any $x\in \mathcal{O}$, $t\in (-\infty,0]$.

 \end{prop}

{\bf Proof} \ \ We can prove along the lines of [Hu2], using the
corresponding result in the manifold case (see [P1, P2] and [CaZ]), and using an orbifold version of Gromoll-Meyer theorem in [Hu2].
  \hfill{$\Box$

\hspace *{0.4cm}

Now we define $\varepsilon$-neck, $\varepsilon$-cap, and strong
$\varepsilon$-neck. As in Definition 2.15 in [KL2], we do not
require the map in the definition of $\varepsilon$-closeness of two
pointed orbifolds to be precisely basepoint-preserving.  Given a
Riemannian 3-orbifold $(\mathcal{O},g)$, an open subset $U$, and a
point $x_0\in U$. $U$ is an $\varepsilon$-neck centered at $x_0$ if
there is a diffeomorphism $\psi: (\mathbb{S}^2//\Gamma)
  \times \mathbb{I} \rightarrow U$  such that the pulled back metric
  $\psi^*g$, scaling with some factor $Q$, is $\varepsilon$-close (in
  $C^{[\varepsilon^{-1}]}$ topology) to the standard metric on $(\mathbb{S}^2//\Gamma)
  \times \mathbb{I}$ with scalar curvature 1 and
  $\mathbb{I}=(-\varepsilon^{-1},\varepsilon^{-1})$, and the distance $d(x_0,  |\psi|(\mathbb{S}^2//\Gamma \times \{0\}))< \varepsilon/ \sqrt{Q}$.
   (Here $\Gamma$ is a finite subgroup of isometries of $\mathbb{S}^2$.)
    An open subset
  $U$ is an $\varepsilon$-cap centered at $x_0$ if  $U$ is diffeomorphic to $\mathbb{R}^3//\Gamma$ or  $\mathbb{S}^2//\Gamma'
\times_{Z_2} \mathbb{R}$, and there is an open set $V$ with
compact closure such that $x_0 \in V \subset \overline{V}\subset U$,
and $U\setminus \overline{V}$ is an $\varepsilon$-neck.
              Given a 3-dimensional  Ricci flow  $(\mathcal{O}, g(t))$, an open subset
$U$, and a point $x_0\in U$.
   $U$ is a strong
  $\varepsilon$-neck  centered at $(x_0,t_0)$ if    there
  is a diffeomorphim $\psi: (\mathbb{S}^2//\Gamma)
  \times \mathbb{I} \rightarrow U$ such that, the pulled back
  solution $\psi^*g(\cdot,\cdot)$ on the parabolic region $\{(x,t)| x \in U, t\in [t_0-Q^{-1},t_0]\}$ (for some $Q>0$),
  parabolically rescaled with factor $Q$, is $\varepsilon$-close (in $C^{[\varepsilon^{-1}]}$
  topology) to the subset $(\mathbb{S}^2//\Gamma
  \times \mathbb{I})\times [-1,0]$ of the evolving round cylinder $\mathbb{S}^2//\Gamma
  \times \mathbb{R}$, with scalar curvature 1 and length 2$\varepsilon^{-1}$ to $\mathbb{I}$ at time zero, and
  the distance at time $t_0$, $d_{t_0}(x_0,  |\psi|(\mathbb{S}^2//\Gamma \times \{0\}))<\varepsilon/ \sqrt{Q}$.

\vspace*{0.4cm} The following proposition is analogous to
Proposition 2.5 in [Hu2].

\begin{prop} \label{prop 2.3} For every $\varepsilon>0$ there exist constants $C_1=C_1(\varepsilon)$ and $C_2=C_2(\varepsilon)$,
such that for every 3-dimensional ancient $\kappa$-solution
$(\mathcal{O},g(\cdot))$  containing no bad 2-suborbifolds, for each space-time point $(x,t)$, there
is a radius $r$, $\frac{1}{C_1}(R(x,t))^{-\frac{1}{2}}< r<
C_1(R(x,t))^{-\frac{1}{2}}$, and an open neighborhood $B$,
$\overline{B_t(x,r)}\subset B \subset B_t(x,2r)$, which falls into
one of the following categories:

(a) $B$ is an strong $\varepsilon$-neck centered at $(x,t)$,

(b) $B$ is an $\varepsilon$-cap,

(c) $B$ is diffeomorphic to a spherical orbifold
$\mathbb{S}^3//\Gamma$  (for a finite subgroup $\Gamma$ of $O(4)$).

Moreover, the scalar curvature in $B$ in cases (a) and (b) at time
$t$ is between $C_2^{-1}R(x,t)$ and $C_2R(x,t)$.
 \end{prop}

{\bf Proof} \ \ The proof is similar to that of Proposition 2.5 in [Hu2],  using
(an extension of) Hamilton's canonical parametrization for a neck region ([H97]).
 \hfill{$\Box$}

\hspace *{0.4cm}

One can easily establish 3-dimensional analogues of Proposition 3.2
and Lemma 3.4  in [Hu2]  on overlapping $\varepsilon$-necks.  We
also  have the following

\begin{prop} \label{prop 2.4}\ \ Let $\varepsilon$ be sufficiently small.
Let $(\mathcal{O},g)$ be a complete, connected 3-orbifold. If each
point of $\mathcal{O}$ is the center of an $\varepsilon$-neck or an
$\varepsilon$-cap, then $\mathcal{O}$ is diffeomorphic to a
spherical orbifold, or a neck, or a cap, or an orbifold connected
sum of at most two spherical 3-orbifolds.
\end{prop}

{\bf Proof}. The proof is along the lines of that of
 Proposition 3.3 in [Hu2].  \hfill{$\Box$}

\section{Ricci flow with surgery on 3-orbifolds}

For the notion of evolving Riemannian orbifold $\{(\mathcal{O}(t), g(t))_{t \in I }\}$ see [Hu2] (also see
[BBM] for the manifold case).

As in [BBM], a time $t\in I$ is regular if $t$ has a neighborhood in
$I$ where $\mathcal{O}(\cdot)$ is constant and $g(\cdot)$ is
$C^1$-smooth. Otherwise it is singular.  We also denote by $f_{max}$
and $f_{min}$ the supremum and infimum of a function $f$,
respectively, as in [BBM].

\hspace *{0.4cm}

  {\bf Definition  }(Compare [BBM], [Hu2])\ \  A piecewise $C^1$-smooth
  evolving Riemannian 3-orbifold $\{(\mathcal{O}(t), g(t))\}_{t \in I }$ is a
  surgical solution of the Ricci flow if it has the following
  properties.

  i. The equation $\frac{\partial g}{\partial t}=-2Ric$ is satisfied
  at all regular times;

  ii.  For each singular time $t_0$ one has $R_{min}(g_+(t_0))\geq
  R_{min}(g(t_0))$;

  iii. For each singular time $t_0$ there is a locally finite collection
  $\mathcal{S}$ of disjoint embedded $\mathbb{S}^2//\Gamma$'s in $\mathcal{O}(t_0)$
  (where $\Gamma$'s are finite subgroups of $O(3)$), and an orbifold $\mathcal{O}'$ such that

  (a) $\mathcal{O}'$ is obtained from $\mathcal{O}(t_0)\setminus \mathcal{S}$ by
  gluing back $\bar{\mathbb{B}}^3//\Gamma$'s;

 (b) $\mathcal{O}_+(t_0)$ is a union of some connected components of $\mathcal{O}'$ and
 $g_+(t_0)=g(t_0)$ on $\mathcal{O}_+(t_0)\cap \mathcal{O}(t_0)$;

(c) Each component of $\mathcal{O}'\setminus \mathcal{O}_+(t_0)$ is
diffeomorphic to a spherical 3-orbifold, or a neck, or a cap, or  an
orbifold connected sum of at most two spherical 3-orbifolds.

\hspace *{0.4cm}

Now we introduce the notion of canonical neighborhood  following
[P2].

\hspace *{0.4cm}

{\bf Definition}  Let $\varepsilon$ and $C$ be positive constants. A
point $(x,t)$ in a surgical solution to the Ricci flow is said to
have an $(\varepsilon,C)$-canonical neighborhood if there is   an open
neighborhood $U$, $\overline{B_t(x,\sigma)} \subset U\subset
B_t(x,2\sigma)$ with $C^{-1}R(x,t)^{-\frac{1}{2}}<\sigma
<CR(x,t)^{-\frac{1}{2}}$, which falls into one of the following four
types:

(a) $U$ is a strong $\varepsilon$-neck with center $(x,t)$,

(b) $U$ is an $\varepsilon$-cap with center $x$ for $g(t)$,

(c) at time $t$, $U$ is diffeomorphic to a closed spherical orbifold
$\mathbb{S}^3//\Gamma$,

\noindent and if moreover, the scalar curvature in $U$ at time $t$  satisfies
the derivative estimates
\begin{equation*}
|\nabla R|< C R^{\frac{3}{2}} \hspace*{8mm} and \hspace*{8mm}
|\frac{\partial R}{\partial t}|< C R^2,
\end{equation*}
and,  for cases (a) and (b), the scalar curvature in $U$ at time $t$
is between $C^{-1}R(x,t)$ and $CR(x,t)$,  and for case (c), the
curvature operator of $U$ is positive, and the infimal sectional
curvature of $U$ is greater than $C^{-1}R(x,t)$.

\hspace *{0.4cm}

 We choose constants $\varepsilon_0>0$ and
$C_0$ similarly as in [Hu2].

Now we consider some a priori assumptions, which consist of the
pinching assumption and the canonical neighborhood assumption.

Following [MT] and [BBB$^+$],  a 3-dimensional surgical solution
$(\mathcal{O}(t), g(t))$ to the Ricci flow
  has curvature pinched toward positive at time $t$ if for any $x \in \mathcal{O}(t)$ there
  holds
\begin{equation*}\begin{aligned}
& R(x,t)\geq -\frac{6}{4t+1}, \\
&  R(x,t)\geq 2(-\nu(x,t))(ln (-\nu(x,t))+ln(1+t)-3) \hspace{4mm}  when \hspace{4mm}  \nu(x,t) < 0, \\
\end{aligned}\end{equation*}
where $\nu$ is the least eigenvalue of the curvature operator.

\hspace *{0.4cm}

{\bf Pinching assumption}:  A  3-dimensional surgical solution
$(\mathcal{O}(t), g(t))$ to the Ricci flow
  satisfies the pinching assumption if it has  curvature pinched toward positive  at all space-time points.

\hspace *{0.4cm}

{\bf Canonical neighborhood assumption}:  Let $\varepsilon_0$ and
$C_0$ be given as above.
 Let $r: [0,+\infty)\rightarrow (0,+\infty)$ be a non-increasing function. An evolving Riemannian 3-orbifold $\{(\mathcal{O}(t), g(t))\}_{t \in I}$
 satisfies the canonical neighborhood assumption  $(CN)_r$ if  any  space-time point $(x,t)$ with  $R(x,t)\geq
r^{-2}(t)$ has an  $(\varepsilon_0,C_0)$-canonical neighborhood.

\hspace *{0.4cm}

Bounded curvature at bounded distance is one of the key ideas in
Perelman [P1], [P2]; compare [MT, Theorem 10.2], [BBB$^+$, Theorem
6.1.1] and [BBM, Theorem 6.4]. 4-dimensional versions have appeared
in [CZ2] and [Hu1, Hu2]. The following version is an extension of
that in [BBB$^+$] and [BBM] to the orbifold case.

\begin{prop} \label{prop 3.1}  For each $A, C >0$ and
each $\varepsilon \in (0,2\varepsilon_0]$, there exists $Q=Q(A,
\varepsilon, C)>0$ and $\Lambda=\Lambda(A, \varepsilon, C)>0$ with
the following property. Let $I\subset [0,\infty)$ and
$\{(\mathcal{O}(t), g(t))\}_{t\in I}$ be a 3-dimensional surgical
solution with bounded curvature (at each time) and satisfying the pinching
assumption. Let $(x_0, t_0)$ be a space-time point such that:

1. $(1+t_0)R(x_0, t_0)\geq Q$;

2. For each point $y\in B(x_0, t_0, AR(x_0, t_0)^{-1/2})$, if $R(y,
t_0)\geq 4R(x_0, t_0)$, then $(y, t_0)$ has an $(\varepsilon,
C)$-canonical neighborhood.

\noindent Then for any $y\in B(x_0, t_0, AR(x_0, t_0)^{-1/2})$, we
have
\begin{equation*}
\frac{R(y, t_0)}{R(x_0, t_0)}\leq \Lambda.
\end{equation*}
\end{prop}

{\bf  Proof} \ \ The proof is similar to that of Proposition 4.1 in
[Hu2].   \hfill{$\Box$}

\hspace *{0.4cm}

The following proposition is analogous to Proposition 4.2 in [Hu2];
compare [BBM, Theorem 6.5] and [BBB$^+$, Theorem 6.2.1].

\begin{prop} \label{prop 3.2}\ \
 For any  $r$, $\delta>0$,
there exist    $h \in (0, \delta r)$ and $D> 10$, such that if
$(\mathcal{O}(\cdot),g(\cdot))$ is a complete 3-dimensional surgical
solution with bounded curvature, defined on a time interval $[a,b]$
and satisfying the pinching assumption and the canonical
neighborhood assumption $(CN)_r$, then the following holds:

 \noindent Let $t \in [a,b]$ and  $x,y, z \in \mathcal{O}(t)$ such that $R(x,t) \leq
2/r^2$,  $R(y,t)=h^{-2}$ and $R(z,t)\geq D/h^2$. Assume there is a
curve $\gamma$ in $\mathcal{O}(t)$ connecting $x$ to $z$ via $y$, such that
each point of $\gamma$ with scalar curvature in $[2C_0r^{-2},
C_0^{-1}Dh^{-2}]$ is the center of an $\varepsilon_0$-neck.  Then
$(y,t)$ is the center of a strong $\delta$-neck.
\end{prop}

{\bf Proof}\ \ The proof is along the lines of that of Proposition
4.2 in [Hu2], using Proposition 3.1.  \hfill{$\Box$}

\hspace *{0.4cm}

The metric surgery on a $\delta$-neck is similar as in [H97] and
[P2]; for the orbifold case see also  the 4-dimensional analog in
[Hu2]. Usually we will be given two non-increasing step functions
$r, \delta: [0,+\infty)\rightarrow (0, +\infty)$ as surgery
parameters. Let $h(r,\delta), D(r,\delta)$ be the associated
parameter  as determined in Proposition 3.2, ($h$ is also called the
surgery scale,) and let $ \Theta:=2Dh^{-2}$ be the curvature
threshold for the surgery process ( as in [BBB$^+$], [BBM] and [Hu1,
Hu2]), that is, we will do surgery when $R_{max}(t)$ reaches
$\Theta(t)$. Now we adapt two more definitions from [BBM] and [Hu1,
Hu2].

\hspace *{0.4cm}

{\bf Definition} (compare [BBM], [Hu1, Hu2] )\ \ Given an interval
$I\subset [0,+\infty)$, fix surgery parameter $r$, $\delta:
I\rightarrow (0,+\infty)$ (two non-increasing functions) and let
$h$, $D$, $\Theta=2Dh^{-2}$ be the associated cutoff parameters. Let
$(\mathcal{O}(t),g(t))$ ($t \in I$) be an evoving Riemannian
3-orbifolds. Let $t_0 \in I$ and $(\mathcal{O}_+,g_+)$ be a
(possibly empty) Riemmanian 3-orbifolds. We say that
$(\mathcal{O}_+,g_+)$ is obtained from
$(\mathcal{O}(\cdot),g(\cdot))$ by $(r,\delta)$-surgery at time
$t_0$ if

i. $R_{max}(g(t_0))=\Theta(t_0)$, and there is a locally finite
collection
  $\mathcal{S}$ of disjoint embedded $\mathbb{S}^2//\Gamma$'s in $\mathcal{O}(t_0)$ which are in the middle  of
  strong $\delta(t_0)$-necks with radius equal to the surgery scale $h(t_0)$, such that
  $\mathcal{O}_+$ is obtained from $\mathcal{O}(t_0)$ by doing
   surgery along these necks ( where $\Gamma$'s are  finite subgroups of $O(3)$
 ), and removing the components each of which is
  diffeomorphic to

  (a)  a spherical 3-orbifold, and either has sectional curvature bounded below by $C_0^{-1}/100$  or
  is
covered by $\varepsilon_0$-necks and $\varepsilon_0$-caps, or

  (b)  a neck, and is  covered by
$\varepsilon_0$-necks, or

 (c) a cap,  and  is covered by $\varepsilon_0$-necks and an $\varepsilon_0$-cap, or

 (d) an orbifold connected sum of at most two spherical 3-orbifolds, and is covered by
$\varepsilon_0$-necks and $\varepsilon_0$-caps.

ii. If $\mathcal{O}_+\neq \emptyset$, then $R_{max}(g_+)\leq
\Theta(t_0)/2$.

\hspace *{0.4cm}

{\bf Definition} (cf. [BBM] and [Hu1, Hu2])\ \  A surgical solution
$(\mathcal{O}(\cdot),g(\cdot))$ defined on some time interval
$I\subset [0,+\infty)$ is an $(r,\delta)$-surgical solution  if it
has the  following properties:

i.  It satisfies the pinching assumption, and $R(x,t) \leq \Theta (t)$ for all $(x,t)$;

ii. At each singular time $t_0\in I$,
$(\mathcal{O}_+(t_0),g_+(t_0))$ is obtained from
$(\mathcal{O}(\cdot),g(\cdot))$ by $(r,\delta)$-surgery at time
$t_0$;

iii. Condition $(CN)_r$ holds.

\hspace *{0.4cm}

 Recall that in our 3-dimensional case, $g(\cdot)$ is
$\kappa$-noncollapsed (for some $\kappa>0$) on the scale $r$ at time
$t$ if at any point $x$, whenever $|Rm|\leq r^{-2} \hspace{2mm} on
\hspace{2mm} P(x, t, r, -r^2) \hspace{2mm}$ we have $ \hspace{2mm}
vol B(x,t, r)\geq \kappa r^3$.
 Let $\kappa: I \rightarrow (0, +\infty)$ be a function.
We say $\{(\mathcal{O}(t), g(t))\}_{t \in I}$ has property $(NC)_\kappa$ if it
is $\kappa (t)$-noncollapsed  on all scales $\leq 1$ at any time
$t\in I$. An $(r,\delta)$-surgical solution which also satisfies
condition $(NC)_\kappa$  is called an $(r,\delta,\kappa)$-surgical
solution.

\hspace *{0.4cm}

Now we establish the existence of Ricci flow with surgery with
initial data a complete 3-orbifold with bounded geometry  and containing no bad 2-suborbifolds. (Note that
by working equivariantly as in [DL], the  conclusion in the
following theorem should also hold true when the initial data has a
covering (instead of itself) which has bounded geometry; but we will
not use this more general version here.)

\begin{thm} \label{thm 3.3}\ \ Given  $v_0>0$, there are surgery parameter
sequences
\begin{equation*}
\mathbf{K}=\{\kappa_i\}_{i=1}^\infty, \hspace{2mm}
\Delta=\{\delta_i\}_{i=1}^\infty, \hspace{2mm}
\mathbf{r}=\{r_i\}_{i=1}^\infty
\end{equation*}
such that the following holds. Let $r(t)=r_i$ and
$\delta(t)=\delta_i$  and $\kappa(t)=\kappa_i$ on $[(i-1)2^{-5},
i\cdot2^{-5})$, $i=1, 2, \cdot\cdot\cdot$.  Let $(\mathcal{O},g_0)$
be  a complete 3-orbifold with  $|Rm|\leq 1$, with vol
$B(x,1)\geq v_0$ at any point $x$, and containing no bad 2-suborbifolds. Then there exists an $(r,\delta,
\kappa)$-surgical solution defined on time interval $[0,\infty)$
with initial data $(\mathcal{O},g_0)$.
\end{thm}

{\bf Proof} The proof is similar to that of Theorem 5.5 in [Hu2].   \hfill{$\Box$}

\hspace *{0.4cm}

Note that if we assume the initial metric has uniformly positive
scalar curvature in addition, the surgical solution constructed in
Theorem 3.3 will become extinct in finite time. Then Theorem 1.1
follows easily, using our procedure of surgery, Proposition 2.4 and
the observations before Proposition 2.2.

\hspace *{0.4cm}

{\bf Acknowledgements}  I'm partially supported by NSFC no.11171025.

\hspace *{0.4cm}

\bibliographystyle{amsplain}

{\bf Reference}

\hspace *{0.4cm}

\bibliography{1}[A] T. Aubin, Some nonlinear problems in Riemannian
geometry, Springer, 1998.

\bibliography{2}[BBB$^+$] L. Bessi$\grave{e}$res, G. Besson, M. Boileau, S.
Maillot and J. Porti, Geometrisation of 3-manifolds, Europ. Math.
Soc. 2010.

\bibliography{3}[BBM] L. Bessi$\grave{e}$res, G. Besson and S. Maillot, Ricci flow on
 open 3-manifolds and positive scalar curvature,  Geometry and Topology  15 (2011), 927-975.

\bibliography{4}[BB]  J. Borzellino,  V. Brunsden,  Elementary orbifold differential
topology, arXiv:1205.1156

\bibliography{5}[CaZ] H.-D. Cao, X.-P. Zhu, A complete proof of the
Poincar$\acute{e}$ and geometrization conjectures- application of
the Hamilton-Perelman theory of the Ricci flow, Asian J. Math. 10
(2006), 165-492.

\bibliography{6}[CTZ] B.-L. Chen, S.-H. Tang and X.-P. Zhu,  Complete classification of
compact four-manifolds with positive isotropic curvature,  J. Diff.
Geom. 91 (2012), 41-80.

\bibliography{7}[CZ] B.-L. Chen, X.-P. Zhu, Ricci flow with surgery
on  four-manifolds with positive isotropic curvature, J. Diff. Geom.
74 (2006), 177-264.

\bibliography{8}[C]  G. Chen, Calculus on orbifolds, Journal of Sichuan University (Natural Science Edition) 41 (2004), no.5, 931-939. (In Chinese)

\bibliography{9}[DL] J. Dinkelbach, B. Leeb, Equivariant Ricci flow with surgery and applications to finite group
actions on geometric 3-manifolds, Geom. Topol. 13 (2009), no.2,
1129¨C1173.

\bibliography{10}[F] C. Farsi, Orbifold spectral theory, Rocky
Mountain J. Math. 31 (2001), no.1, 215-235.

\bibliography{11}[GL] M. Gromov, H. B. Lawson, Jr, The classification of simply connected
manifolds of positive scalar curvature, Ann. of Math. 111 (1980),
423-434.

\bibliography{12}[H82] R. Hamilton, Three-manifolds with positive
Ricci curvature, J. Diff. Geom. 17(1982), 255-306.

\bibliography{13}[H97] R. Hamilton, Four-manifolds with positive isotropic curvature, Comm. Anal. Geom. 5 (1997), 1-92.

\bibliography{14}[He] E. Hebey,
Nonlinear analysis on manifolds: Sobolev spaces and inequalities,
Courant Institute of Mathematical Sciences, 1999

\bibliography{15}[Hu1] H. Huang, Ricci flow on open 4-manifolds with positive isotropic
curvature, arXiv:1108.2918, to appear in J. Geom. Anal.; published
online first,   DOI 10.1007/s12220-011-9284-y

\bibliography{16}[Hu2] H. Huang,  Complete 4-manifolds with uniformly positive isotropic
curvature, arXiv:0912.5405v12

\bibliography{17}[KL1] B. Kleiner, J. Lott, Notes on Perelman's
papers, Geom. Topol. 12 (2008), 2587-2855;  arXiv:0605667v4.

\bibliography{18}[KL2]  B. Kleiner, J. Lott, Geometrization of three-dimensional orbifolds via Ricci flow, arXiv:1101.3733v2.

\bibliography{19}[Ko]  B. Kotschwar, An energy approach to the problem of uniqueness for the Ricci
flow, arXiv:1206.3225.

\bibliography{20}[MT] J. Morgan, G. Tian, Ricci flow and the
Poincar$\acute{e}$ conjecture, Clay Mathematics Monographs 3, Amer.
Math. Soc., 2007.

\bibliography{21}[N] Y. Nakagawa, An isoperimetric inequality for
orbifolds, Osaka J. Math. 30 (1993), 733-739.

\bibliography{22}[P1] G. Perelman, The entropy formula for the Ricci flow and its geometric applications,
arXiv:math.DG/0211159.

\bibliography{23}[P2] G. Perelman, Ricci flow with surgery on three-manifolds, arXiv:math.DG/0303109.

\bibliography{24}[P3]  G. Perelman, Finite extinction time for the
solutions to the Ricci flow on certain three-manifolds,
arXiv:math.DG/0307245.

\bibliography{25}[SY] R. Schoen, S.-T. Yau, On the structure of manifolds with positive scalar
curvature, Manuscripta Math. 28 (1979),  no. 1-3, 159-183.

\bibliography{26}[S] W.-X. Shi, Deforming the metric on complete
Riemannian manifolds, J. Diff. Geom. 30 (1989), 223-301.

\bibliography{27}[T] W. Thurston, The geometry and topology of 3-manifolds, Princeton lecture notes (1979).

\bibliography{28} [Z]  Q. S. Zhang, Sobolev inequalities, heat
kernels under Ricci flow, and the Poincar$\acute{e}$ conjecture, CRC
Press 2011.

\vspace *{0.4cm}

School of Mathematical Sciences, Key Laboratory of Mathematics and Complex Systems,

Beijing Normal University, Beijing 100875, P.R. China

 E-mail address: hhuang@bnu.edu.cn

\end{document}